# *Une initiation au concept de l'infini : Genèse et développements dus à Cantor*

Michel Adès*, David Guillemette*, Serge B. Provost**
*UQAM, **Université Western Ontario

## Résumé

On s'intéresse dans cet article à explorer la notion de l'infini en étudiant la contribution de Cantor à ce domaine. On fait un bref rappel historique sur la théorie des ensembles. Comme exemple de l'infini, on considère le célèbre hôtel de Hilbert. On fait une construction graphique pour démontrer la dénombrabilité des rationnels. On approfondit le procédé diagonal de Cantor qui démontre que l'infini des réels est plus grand que celui des entiers. Il y a donc une hiérarchie de l'infini, ce qui nous permet d'élaborer sur l'hypothèse du continu qui demeure toujours un problème non résolu en mathématiques.

## Adresses postales

Michel Adès / David Guillemette            Serge Provost
Université du Québec à Montréal            The University of Western Ontario
Département de mathématiques               Département de statistique et d'actuariat
Local PK-5151                              Western Science Centre, Room 262
201, Avenue du Président-Kennedy           London (Ontario)
Montréal (Québec) H2X 3Y7                  Canada N6A 5B7

## Auteur correspondant
Michel Adès

## Courriels
ades.michel@uqam.ca / guillemette.david@uqam.ca / sp@uwo.ca

# *Une initiation au concept de l'infini : Genèse et développements dus à Cantor*

Michel Adès*, David Guillemette*, Serge B. Provost**
*UQAM, **Université Western Ontario

L'infini est un sujet qui passionne l'être humain depuis l'Antiquité. Par exemple, est-ce que notre univers est fini ou infini ? Est-il en expansion perpétuelle vers un inaccessible infini?

## 1- Introduction

Georg Cantor, mathématicien allemand d'origine russe, naquit le 3 mars 1845 à St-Pétersbourg en Russie et décéda le 6 janvier 1918 à Halle, en Allemagne. Il est reconnu comme étant le fondateur de la théorie des ensembles. Cantor est considéré comme le mathématicien et le législateur de l'infini, car il développa une arithmétique de l'infini. Cantor aurait identifié « l'infini absolu » à Dieu, comme étant l'infini d'une classe propre comme celle de tous les cardinaux.

L'infini est un château-fort qui a dévoilé quelques-uns de ses mystères à Cantor, lequel obtint des résultats prodigieux, ce qui lui a causé beaucoup d'inconforts et attiré les sarcasmes de la part des grands mathématiciens de son époque, tels que Leopold Kronecker (1823-1891) et Henri Poincaré (1854-1912).

Sur le plan étymologique, le mot infini est d'origine latine « infinitus » qui signifie sans limites. Le symbole de l'infini représenté par «∞» a été introduit pour la première fois en 1655 par le mathématicien anglais John Wallis (1616-1703) dans son ouvrage « *De sectionibus conicis* » : Mathématiques sur les sections coniques. Notons que ce symbole a été popularisé vers 1713, grâce à Jacques Bernoulli qui est l'inventeur, entre autres, de la fameuse lemniscate qui signifie ruban.

## 2- Brève perspective historique

Entre le 17$^e$ et 19$^e$ siècle, les mathématiques vont connaître un développement fulgurant en Europe. La création de nombreuses académies et institutions scientifiques, ainsi que l'apparition de pôles d'enseignement importants, comme les écoles polytechniques, vont contribuer à cet essor très important, voire trop important ! Car, pour plusieurs, ces développements se sont parfois faits sans véritables fondations solides. Fondations sur lesquelles ces développements peuvent s'asseoir et au regard desquelles ils peuvent être validés. On pense tout de suite en exemple au développement du calcul différentiel et intégral, qui a été sans doute l'une des plus grandes avancées mathématiques du 17$^e$ siècle. Les notions fondamentales de continuité, de dérivabilité et même de fonctions doivent encore à l'époque de Cantor être clarifiées.

C'est une époque où apparaissent aussi de nombreux objets nouveaux et notions originales qui sont difficiles à appréhender, à classer et à articuler aux différents champs mathématiques. On peut donner en exemple la notion de groupe, de vecteur, de matrice ou encore de géométries non-euclidiennes. Apparaissent aussi des objets énigmatiques et paradoxaux qui appellent à des questions de fondements comme les fonctions continues et pourtant dérivables en aucun point. L'apparition et l'utilisation parfois féconde de plusieurs notions et objets originaux, susciteront la controverse et inciteront les mathématiciennes et mathématiciens de l'époque à une réflexion orientée davantage en termes de fondements.

### 3- **Théorie des ensembles : Genèse et développement**

En 1872, Georg Cantor rencontre Richard Dedekind en Suisse. C'est le début d'une amitié solide et d'une longue collaboration. Suivra une intense correspondance qui a fortement renseigné les historiens des mathématiques quant au début de la théorie des ensembles. Cette théorie cherche à fournir un langage et des notions communes à divers champs mathématiques dans une visée unifiante et englobante. Autrement dit, une théorie capable de fonder et de donner du sens à l'ensemble des développements mathématiques de cette époque.

L'idée simple d'ensemble en tant que « collection d'objets que nous pouvons réunir par la pensée » est au cœur de cette théorie qui marque les mathématiques.

Dès le départ, les ensembles de nombres et leurs relations vont intéresser Cantor et Dedekind. Cantor introduit la notion d'ensembles dits *de même puissance*. Il fait l'énoncé suivant : « Si deux ensembles bien définis M et N se laissent coordonner l'un à l'autre, élément par élément, de façon univoque et complète, je me sers de l'expression… qu'ils sont de la même puissance ou qu'ils sont équivalents[1] ».

En fait, deux ensembles sont équivalents ou équipotents lorsqu'il existe une bijection entre eux. On dira alors qu'ils ont la même cardinalité ou le même nombre d'éléments.

L'idée de Cantor est d'établir une relation biunivoque entre les éléments de deux ensembles, autrement dit, de ranger ou d'associer par une règle précise chaque élément d'un premier ensemble avec un unique élément d'un second ensemble, pour en arriver à la conclusion que ces ensembles sont constitués d'un même nombre d'éléments. Par exemple, imaginons qu'un grand nombre de personnes assistent à un concert. Pour connaître le nombre de personnes dans l'assistance, on n'est pas obligé de les compter une à une. Si nous connaissons la capacité de la salle, on peut regarder le nombre de chaises inoccupées ou le nombre de personnes sans chaise. Dans le cas où tout le monde est assis, c'est-à-dire dans le cas où chaque personne sans exception est associée à une unique chaise et qu'il n'y a pas de chaise vide, on peut conclure que le nombre de personnes correspond à la capacité de la salle.

---

1. *cité dans Dahan-Dalmedico* & *Peiffer, 1986, p. 239.*

On a, en fait, réussi à associer de façon biunivoque chaque élément de l'ensemble contenant des personnes à chaque élément de l'ensemble contenant les chaises[2].

Plus précisément, on dira avoir établi une relation bijective entre les deux ensembles. Une relation bijective est une relation à la fois injective et surjective.

La relation est ici injective, car tous les éléments de l'ensemble de départ sont associés à un et un seul élément de l'ensemble d'arrivée (chaque personne est assise et est assise sur une unique chaise). Elle est aussi surjective, car tous les éléments de l'ensemble d'arrivée sont associés à au moins un élément de l'ensemble de départ (chaque chaise est occupée par au moins une personne).

Ainsi, pour établir que deux ensembles ont le même nombre d'éléments, il suffit d'établir une bijection entre eux. Le plus étonnant, et souvent contre-intuitif, est qu'il en ira de même pour les ensembles contenant un nombre infini d'éléments. Comme nous l'aborderons plus loin en détail, Cantor montrera en premier lieu que l'on peut ranger l'ensemble des rationnels $\mathbb{Q}$ en une suite simple et indexée par les entiers naturels $\mathbb{N}$. Ce qui veut dire que l'on peut établir une relation bijective entre $\mathbb{N}$ et $\mathbb{Q}$, et donc qu'il y aurait « autant » de nombres naturels que de nombres rationnels, quoiqu'en nombre infini… On dira aujourd'hui que $\mathbb{N}$ et $\mathbb{Q}$ ont la même *cardinalité*, quelque chose comme « la même taille ».

Les ensembles ayant la même cardinalité que $\mathbb{N}$ sont dits *dénombrables*, c'est-à-dire qu'ils peuvent être « comptés » à l'aide des nombres naturels. Dans cette foulée, Dedekind montrera ensuite que l'ensemble des nombres algébriques est aussi dénombrable. L'affaire était lancée… Quels ensembles sont-ils dénombrables? Existe- t il des ensembles indénombrables?

En fait, un nombre θ est dit un nombre algébrique[3], lorsqu'il satisfait à une équation $\theta^n + a_1\theta^{n-1} + a_2\theta^{n-2} + \ldots + a_{n-1}\theta + a_n = 0$, de degré fini $n$ et à coefficients rationnels $a_1, a_2, \ldots, a_{n-1}, a_n$.

Comme exemples d'un nombre algébrique, on considère les cas suivants : a) Tout nombre rationnel p/q est algébrique car le quotient de deux entiers est solution de $qx - p = 0$; b) Le nombre imaginaire ou complexe i est algébrique, car il est solution de l'équation $x^2 + 1 = 0$; c) Racine carrée de 2, qui est solution de l'équation $x^2 - 2 = 0$, est algébrique et aussi irrationnel; d) Le nombre d'or est également un nombre algébrique défini comme étant l'unique solution positive de l'équation $x^2 - x - 1 = 0$, dont la valeur est 1,6180339887. . .

Maintenant, comment deux ensembles, dont l'un est inclus dans l'autre, peuvent-ils avoir la même taille? La question étonnera les deux mathématiciens eux-mêmes dès le départ. À ce sujet, Cantor écrira à Dedekind : « Je le vois, mais je ne le crois pas! ». On verra bientôt la réponse à ce paradoxe.

2. *Nous remercions M. André Boileau, professeur associé au Département de mathématiques de l'UQAM, à qui nous devons cet exemple.*
3. *R. Dedekind (voir bibliographies).*

## 4- L'hôtel infini de Hilbert

Après ce bref discours sur l'infini, Hilbert nous invite à visiter son hôtel paradoxal avant de poursuivre notre chemin en quête d'infini et de nouveaux savoirs.

Supposons que vous êtes directeur d'un hôtel qui affiche *complet*. Dans la vraie vie, il n'est donc plus possible d'accueillir de nouveaux clients. En revanche, dans le monde des mathématiques, ça ne pose pas de problème si l'hôtel dispose d'un nombre infini de chambres. Cette idée remonte au mathématicien allemand David Hilbert qui a utilisé l'exemple d'un tel hôtel pour démontrer les jeux contre-intuitifs auxquels on peut s'adonner avec le concept de l'infini.

**4.1)** Supposons donc qu'un hôtel dispose d'une infinité de chambres, numérotées 1, 2, 3, …, et que toutes ces chambres sont déjà occupées. Que faire alors si un nouveau client se présentait? Il suffira simplement de demander à la personne de la chambre 1 de déménager dans la chambre 2, à celle de la chambre 2 d'emménager dans la chambre 3, à celle de la chambre 3 de se déplacer dans la chambre 4, etc.

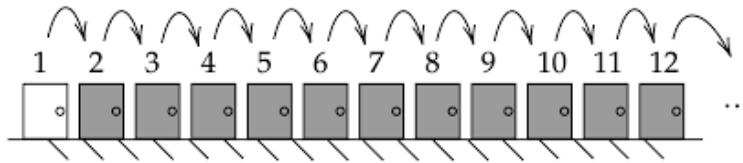

S'il n'y avait qu'un nombre limité de chambres, la personne résidente de la dernière chambre n'aurait nulle part où aller, mais comme il y en a une infinité, toutes et tous auront une chambre où se loger. (On devra cependant demander aux personnes résidentes de se déplacer simultanément, car si elles le font l'une après l'autre, le déménagement peut requérir un temps infini puisqu'une infinité de personnes résidentes doivent changer de chambres!).

**4.2)** En utilisant cette manœuvre astucieuse, on peut en fait accueillir n'importe quel nombre fini de nouvelles personnes invitées. Par exemple, advenant que *p* nouvelles personnes se présentaient à la réception, on demandera alors à chaque personne résidente de l'hôtel de se déplacer vers la chambre dont le numéro est celui de sa chambre plus *p*. Les deux prochaines étapes soulignent une affluence vertigineuse des autocars et leurs passagers.

**4.3)** Supposons maintenant qu'un nombre infini de nouvelles personnes désirent se loger à cet hôtel. Dans ce cas, il suffira de demander à chaque personne résidente d'emménager dans la chambre dont le numéro est le double de celle qu'elle occupe présentement. Une fois cette manœuvre complétée, seules les chambres paires seront occupées, et les chambres impaires pourront alors accueillir un nombre infini de nouvelles personnes.

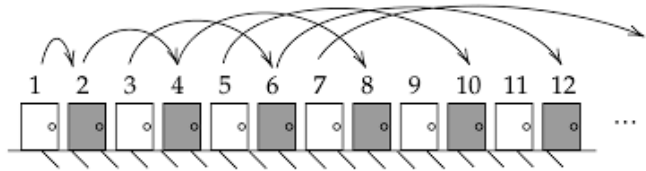

**4.4)** Un peu plus tard, arrive un groupe encore plus important constitué d'une infinité d'autocars, chacun ayant à son bord une infinité de passagers.

Le maître d'hôtel demande alors à la personne résidente de la chambre j ≥ 0 de se rendre à la chambre 2j+1, toutes les chambres possédant un numéro pair devenant ainsi libres; de plus, il donne la consigne qui suit aux responsables des autocars : *le passager **j** de l'autocar **i** ≥ **0*** doit occuper la chambre $(2j + 1) \, 2^{(i+1)}$. Chaque passager aura ainsi sa propre chambre. En effet :

**4.4.1)** Pour i=0 et j ≥ 0, les chambres sont : 2, 6, 10, 14, 18, 22, 26, 30, 34, 38, . . .

**4.4.2)** Pour i= 1 et j ≥ 0, les chambres sont : 4, 12, 20, 28, 36, 44, 52, 60, 68, . . .

**4.4.3)** Pour i = 2 et j ≥ 0, les chambres sont : 8, 24, 40, 56, 72, 88, 104, 120, 136, …

Pour les cas 4.4.1, 4.4.2 et 4.4.3, les chambres sont respectivement le double, le quadruple et l'octuple d'un nombre impair, et les prochaines sont 16 fois, 32 fois un nombre impair, etc. En général, les chambres sont $2^{i+1}$ fois un nombre impair pour i ≥ 0.

Notons que tous les voyageurs auront leur propre chambre car deux voyageurs distincts ne pourront se retrouver dans la même chambre. À ce sujet, la décomposition en facteurs premiers consiste à écrire un entier naturel sous la forme d'un produit de nombres premiers. En fait, le théorème fondamental de l'arithmétique affirme que tout entier strictement positif possède une unique décomposition en facteurs premiers. Ainsi, pour *j* différent de *r* et *i* différent de *s*, l'expression $(2j + 1) \, 2^{(i+1)}$ est différente de $(2r + 1) \, 2^{(s+1)}$, car $(i + 1)$ est l'exposant de 2 dans la décomposition en facteurs premiers qui est unique.

## 5- Les nombres naturels, entiers et rationnels

Nous rappelons maintenant quelques notions élémentaires sur les nombres qui seront très utiles dans cet article.

Soit $\mathbb{N}$ = {0, 1, 2, 3, 4, 5, 6, 7, . . .} l'ensemble des nombres naturels ($\mathbb{N}$ étant l'initiale de naturel) dont on sait depuis Euclide (mathématicien grec qui a vécu entre 325-265 avant Jésus-Christ) que le cardinal ou le nombre de ses éléments est infini. En fait, c'est l'infini le plus simple et le plus petit parmi les infinis, que Cantor appela $\aleph_0$, prononcé aleph zéro. Le cardinal de $\mathbb{N}$ est # $\mathbb{N}$ = $\aleph_0$. Notons que $\aleph_0$ est le cardinal de l'ensemble dénombrable des nombres entiers et rationnels. L'indice zéro est utilisé pour désigner le plus petit de tous les aleph (aleph étant la première lettre de l'alphabet hébreu et arabe).

Soit $\mathbb{Z}$ = {. . ., -5, -4, -3, -2, -1, 0, 1, 2, 3, 4, 5, . . .} l'ensemble des nombres entiers relatifs ou simplement entiers, et soit $\mathbb{Q}$ = { $\frac{p}{q}$ | $p \in \mathbb{Z}$, $q \in (\mathbb{Z} \setminus \{0\})$} l'ensemble des rationnels.

De l'ensemble $\mathbb{N}$, prenons par exemple les nombres premiers $\mathbb{P}$ = {2, 3, 5, 7, 11, . . .}, et établissons une relation bijective entre $\mathbb{N}$ et $\mathbb{P}$ :

```
0   1   2   3   4   5   6   7   8   9   10  11 . . .
↕   ↕   ↕   ↕   ↕   ↕   ↕   ↕   ↕   ↕   ↕   ↕ . . .
2   3   5   7   11  13  17  19  23  29  31  37 . . .
```

Il est clair que $\mathbb{P}$ est une partie de l'ensemble infini $\mathbb{N}$ qui est à son tour infini, et pourtant, $\mathbb{P}$ et $\mathbb{N}$ ont la même cardinalité $\aleph_0$, ce qui semble paradoxal !

Ce qui est également plus paradoxal et même stupéfiant est que $\mathbb{Z}$ et $\mathbb{N}$ ont la même cardinalité, bien que $\mathbb{N} \subset \mathbb{Z}$, en effet :

```
0   1   2   3   4   5   6   7   8   9   10 . . .
↕   ↕   ↕   ↕   ↕   ↕   ↕   ↕   ↕   ↕   ↕ . . .
0  -1   1  -2   2  -3   3  -4   4   5  -5 . . .
```

Ainsi, dans ce contexte de dénombrabilité infinie, un sous-ensemble infini a le même cardinal que l'ensemble infini duquel il est extrait.

Les cardinaux infinis n'ont pas la même arithmétique ou le même comportement que les cardinaux finis; en fait, 1 + 1 = 2, alors que $\aleph_0 + \aleph_0 = \aleph_0$, ce qui est étrange, mais pourtant vrai. En effet, $\mathbb{Z}$ est de cardinalité $\aleph_0$ alors qu'il est l'union des ensembles entiers positifs et négatifs dont chacun est de cardinalité $\aleph_0$.

Maintenant pour ce qui est de l'ensemble $\mathbb{Q}$, il existe plusieurs fractions positives ou négatives qui donnent le même nombre, comme par exemple $\frac{24}{6} = \frac{16}{4} = \frac{8}{2} = 4$.

De plus, la plupart des fractions ne sont pas des nombres entiers. Il y a donc plus de fractions que d'entiers, nous avons ainsi $\mathbb{N} \subset \mathbb{Z} \subset \mathbb{Q}$ qui sont d'ailleurs tous dénombrables.

## 6- La dénombrabilité des rationnels

Comme $\mathbb{Q}$ contient l'ensemble $\mathbb{N}$, on peut s'attendre intuitivement à ce que $\mathbb{Q}$ contienne plus d'éléments que $\mathbb{N}$, or, ce n'est pas le cas. En effet, Cantor a démontré que la cardinalité de $\mathbb{Q}$ est la même que celle de $\mathbb{N}$, ce qui semble de prime abord contre-intuitif et même aberrant. En voici la représentation graphique pour les rationnels positifs illustrée sous la forme d'un certain réseau, où on associe un nombre entier à un rationnel :

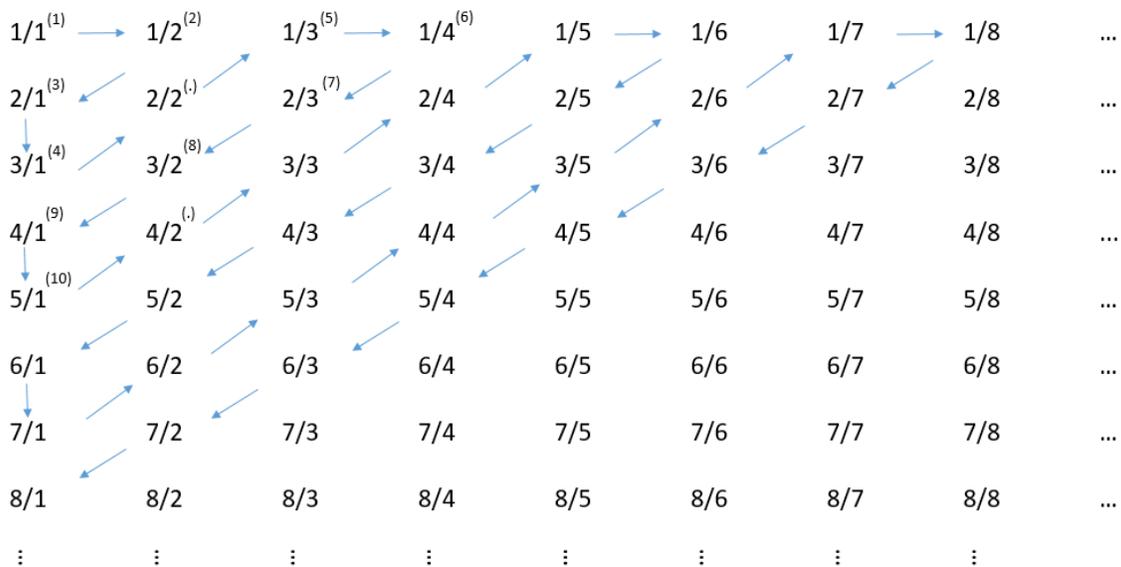

Cette construction graphique contient toutes les fractions ou nombres rationnels. En effet, on observe qu'on a tous les numérateurs et dénominateurs possibles sur chaque colonne et ligne respectivement. Pour repérer un nombre rationnel ou une fraction $\frac{p}{q}$ sur ce graphique, il faut, en suivant les flèches numérotées {1, 2, 3, 4, 5, 6, 7, 8, 9, . . .}, s'acheminer à la $p$-ième ligne et $q$-ième colonne pour atteindre le nombre rationnel voulu. En dénombrant ladite suite et en omettant de mettre un numéro aux rationnels qui sont déjà apparus précédemment sur la liste, on obtient un merveilleux résultat; l'égalité des cardinaux de $\mathbb{Q}$, $\mathbb{N}$ et $\mathbb{Z}$, c'est-à-dire :
$\# \mathbb{Q} = \# \mathbb{N} = \# \mathbb{Z} = \aleph_0$.

## 7- Les nombres réels

Tout d'abord, un nombre irrationnel est un nombre qui ne peut pas s'écrire sous la forme d'une fraction, tels, $\sqrt{2}, \sqrt{3}, \pi, e$. Un nombre réel est un nombre rationnel ou irrationnel, l'ensemble des nombres réels est représenté par $\mathbb{R}$. Notons qu'à partir de n'importe quel irrationnel, on peut former une infinité d'irrationnels, comme $\{\sqrt{2}, 2\sqrt{2}, 3\sqrt{2}, ...\}$, $\{e, 2e, 3e, ...\}$, $\{-\sqrt{5}, -2\sqrt{5}, -3\sqrt{5}, ...\}$. On peut donc construire une infinité d'infinités; ainsi, la cardinalité $\aleph_1$ de $\mathbb{R}$ est plus grande que celles de $\mathbb{N}$ et $\mathbb{Q}$. La lettre $\aleph_1$, prononcé aleph 1, correspond à la notation de Cantor. Précisons que $\aleph_1$ désigne le cardinal de l'ensemble non-dénombrable des nombres réels et des irrationnels.

## 8- Le procédé diagonal et la non-dénombrabilité des réels

Cantor a montré qu'on ne peut pas établir une correspondance « un à un » entre les réels et les entiers. Sa démonstration de la non-dénombrabilité des réels est basée sur le procédé diagonal. Elle consiste à postuler une hypothèse qui va s'avérer fausse, ce qui en fait une preuve par l'absurde. En fait, si *tous* les réels dans l'intervalle (0, 1) étaient en correspondance biunivoque avec les entiers positifs, on pourrait alors tous les énumérer comme ci-contre par exemple.

1 ↔ 0,**7**635789675...
2 ↔ 0,3**5**67000000...
3 ↔ 0,76**8**5643219...
4 ↔ 0,679**4**583210...
5 ↔ 0,2154**3**98607...
6 ↔ 0,11546**2**9703...
7 ↔ 0,874230**1**784...
8 ↔ 0,9567012**9**68...
9 ↔ 0,30976211**5**2...
10 ↔ 0,267451487**9**...
...
...

Dans la configuration listée, tous les réels dans (0,1) doivent y figurer en principe. Alors, l'idée ingénieuse de Cantor en 1874 est de prendre les chiffres en gras sur la **diagonale** et de changer leur valeur en ajoutant 1 s'ils sont dans l'intervalle [0, 8] et de remplacer 9 par 0.

On obtient donc un nouveau chiffre soit **0,8695432060...** qui ne peut pas être dans la liste car il diffère de chacun des réels qui s'y trouvent déjà. En détails, ce nombre **0,869 543 206 0...**, que l'on vient de construire artificiellement ne peut figurer sur la liste des nombres qui s'y trouvent déjà, car son premier chiffre (le **8**) est différent de celui du premier nombre déjà existant sur la liste (0,**7**635789675...). De plus, son deuxième chiffre (le **6**) est différent de celui du deuxième nombre déjà existant sur la liste (0,3**5**67000000...), etc. On peut continuer ce procédé diagonal une infinité de fois. On constate que ce nombre ne peut figurer sur la liste car son nième chiffre est toujours différent du nième chiffre du nième nombre déjà sur la liste.

Ainsi, l'hypothèse de départ, à savoir que tous les réels dans (0, 1), dénombrés et répertoriés les uns en-dessous des autres, est fausse. En conséquence, les réels dans (0, 1) ne sont pas dénombrables, et il y a plus de réels que d'entiers. On a donc $\aleph_0 < \aleph_1$ tel que $\# \mathbb{R} = \aleph_1$.

### 9- L'hypothèse du continu

On vient de montrer que l'infini des rationnels est de cardinalité $\aleph_0$ qui est bien dénombrable alors que l'infini des réels est de cardinalité $\aleph_1$ qui est non-dénombrable, de plus, $\aleph_0 < \aleph_1$.

Cantor a énoncé en 1877 une hypothèse selon laquelle il n'existe pas d'ensemble à la fois plus grand que $\mathbb{N}$ et plus petit que $\mathbb{R}$. C'est l'hypothèse du continu, selon laquelle il n'existe pas d'ensemble infini $\mathbb{E}$ intermédiaire entre $\mathbb{N}$ et $\mathbb{R}$ tel que $\{\aleph_0 < \# \mathbb{E} < \aleph_1\}$. Ladite hypothèse figurait en première position sur la célèbre liste des 23 problèmes de David Hilbert qu'il formula au tout début du vingtième siècle pour la présenter au congrès international des mathématiques qui s'est tenu à Paris en 1900.

Kurt Gödel (1906–1971) a démontré en 1938 que l'hypothèse du continu n'est pas réfutable dans la théorie des ensembles de ZFC[4], généralement utilisée en mathématiques et considérée comme une formalisation adéquate de la théorie des ensembles de Cantor. Ladite hypothèse est donc indépendante des axiomes de ZFC, ou encore indécidable dans le cadre de cette théorie.

Paul Cohen (1934–2007) a démontré en 1963 avec sa propre méthode appelée « forcing » que l'hypothèse du continu est indépendante des axiomes standards de la théorie des ensembles de ZFC. Ainsi, il est impossible de décider s'il existe un ensemble $\mathbb{E}$ entre $\mathbb{N}$ et $\mathbb{R}$ ou pas, c'est l'indécidabilité. On peut donc supposer que l'hypothèse du continu est vraie ou fausse sans obtenir de contradiction dans cette théorie. L'hypothèse du continu reste toujours un sujet actif de recherche en théorie des ensembles. L'ajout de nouveaux axiomes (axiome de détermination, axiomes de grands cardinaux etc.) au corpus de la théorie de ZFC va peut-être aider un jour, grâce aux travaux récents de William Hugh Woodin (1955–), à jeter une lumière spécifique et nouvelle sur la fameuse hypothèse du continu et possiblement la démontrer.

---

4. *Zermelo-Fraenkel-Cantor, Ernst Zermelo (1871-1953) et Abraham Fraenkel (1891-1965).*

## 10-Jean Dieudonné et l'hypothèse du continu

En terminant cette section, on se réfère à cette belle citation de l'éminent mathématicien français Jean Dieudonné à propos de l'hypothèse du continu : « Ce qui nous intéresse beaucoup, en tant que garde-fou, ce sont les preuves d'indécidabilité et d'impossibilité. Des mathématiciens ont passé des années de leur vie à essayer de démontrer l'hypothèse du continu, problème qui les a tourmentés pendant très longtemps. Je me souviens d'avoir entendu dire de mon maître George Pólya, qui le tenait lui-même d'Alexandroff, qu'Alexandroff avait pendant un an travaillé à la démonstration de l'hypothèse du continu et puis qu'il avait arrêté parce qu'il se sentait devenir fou. Il a bien fait. Alors, quand Gödel et Cohen sont venus nous dire qu'il était inutile de nous tracasser les méninges et que jamais nous ne démontrerions ni l'hypothèse du continu ni sa contradiction, nous avons dit : Ouf … ».

## 11- Quelques compléments

**11.1)** Considérons un ensemble F contenant $n$ éléments. En utilisant l'argument de la diagonale, Cantor démontre dans son théorème de 1891, que l'ensemble des parties de l'ensemble F noté par P(F) de cardinal $2^n$, a toujours strictement plus d'éléments que F, même si F est infini. Autrement dit, il n'existe pas de bijection entre F et P(F). La conséquence de ce théorème est l'existence d'une hiérarchie stricte, elle-même infinie, d'ensembles infinis.

**11.2)** L'hypothèse du continu stipule que $c = \aleph_1 = 2^{\aleph_0}$ est l'infini associé à l'ensemble des nombres réels, ce qui signifie que $\aleph_1$ est le nombre de l'ensemble des parties de l'ensemble des entiers. Ce théorème de Cantor traduit le fait qu'il y a des infinis plus grands que d'autres, et que de plus, il n'existe aucune cardinalité infinie qui soit la plus grande.

## 12- Conclusion

Cantor subira les attaques et les sarcasmes de ses contemporains. En effet, Kronecker, qui fut son ancien professeur à Berlin, l'accusa de corrompre la jeunesse et le traita de renégat et de charlatan, alors que Poincaré pensait que les idées de Cantor constituaient une grave maladie des mathématiques, une affection perverse qui serait un jour guérie. En contrepartie, Hilbert a qualifié le travail de Cantor de paradis mathématique « Nul ne doit nous exclure du Paradis que Cantor a créé ». D'ailleurs, Bertrand Russel considérait Cantor comme l'un des plus grands penseurs du XIXe siècle (cité dans Freiberger & Thomas, 2015, p. 251).

## 13- Bibliographie

([https://gdz.sub.unigoettingen.de/id/PPN237853094?tify=%7B%22pages%22%3A%5B19%5D%2C%22view%22%3A%22info%22%2C%22zoom%22%3A0.479%7D](https://gdz.sub.unigoettingen.de/id/PPN237853094?tify=%7B%22pages%22%3A%5B19%5D%2C%22view%22%3A%22info%22%2C%22zoom%22%3A0.479%7D)).

. Dahan-Dalmedico, A. & Peiffer, J., *« Une histoire des mathématiques : routes et dédales ».* Seuil, Paris, 1986.

. Freiberger, M. & Thomas, R., « Dans le secret des nombres ». Dunod, Paris, 2015.

. Wallis, John, « De sectionibus conicis nova methodo expositis tractatus », 1655. ([https://archive.wikiwix.com](https://archive.wikiwix.com))

************************************************

# Biographies : un bref aperçu

Georg Cantor

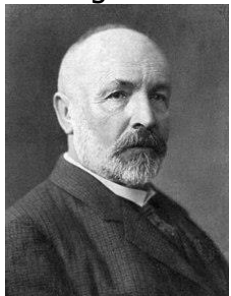

(1845-1918)

Georg Cantor, mathématicien allemand d'origine russe, naquit le 3 mars 1845 à St-Pétersbourg en Russie et décéda le 6 janvier 1918 à Halle, en Allemagne. Il est reconnu comme étant le fondateur de la théorie des ensembles. Cantor est considéré comme le mathématicien et le législateur de l'infini, car il développa une arithmétique de l'infini. Cantor aurait identifié « l'infini absolu » à Dieu, comme étant l'infini d'une classe propre comme celle de tous les cardinaux.

------------------------------------------------

## Richard Dedekind

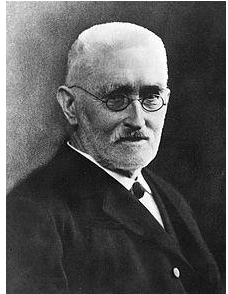

(1831-1916)

Mathématicien allemand, pionnier de l'axiomatisation de l'arithmétique, il a proposé une définition axiomatique de l'ensemble des nombres entiers ainsi qu'une construction rigoureuse des nombres réels à partir des nombres rationnels soit la méthode des « coupures » de Dedekind.

---

## David Hilbert

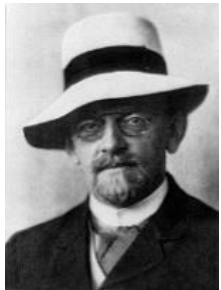

(1862- 1943)

Mathématicien allemand, considéré comme l'un des plus grands mathématiciens du XXe siècle. Il a créé ou développé un large éventail d'idées fondamentales, que ce soit la théorie des invariants, l'axiomatisation de la géométrie ou les fondements de l'analyse fonctionnelle. L'un des exemples les mieux connus de sa position de chef de file est sa présentation, en 1900, de ses fameux 23 problèmes qui ont durablement influencé les recherches mathématiques du XXe siècle. Hilbert a développé une portion significative de l'infrastructure mathématique nécessaire à l'éclosion de la mécanique quantique et de la relativité générale.

---------------------------------------------------------------------------

## Henri Poincaré

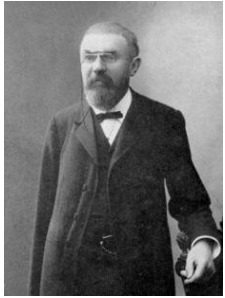

(1854- 1912)

Mathématicien français, dont les avancées sur le problème des trois corps en font un fondateur de l'étude qualitative des systèmes d'équations différentielles et de la théorie du chaos. Précurseur majeur de la théorie de la relativité restreinte et de la théorie des systèmes dynamiques. Poincaré est considéré comme un des derniers grands savants universels, maîtrisant l'ensemble des branches des mathématiques de son époque et certaines branches de la physique.

---------------------------------------------------------------------------

## John Wallis

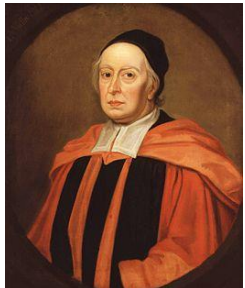

(1616-1703)

Mathématicien anglais dont les travaux sont précurseurs de ceux de Newton. Il a introduit le symbole de l'infini, représenté par ∞, pour la première fois en 1655 dans son ouvrage « De sectionibus conicis » : Mathématiques sur les sections coniques. Ce symbole a été popularisé vers 1713 grâce à Jacques Bernoulli qui est l'auteur, entre autres, de la fameuse lemniscate qui signifie ruban.

----------------------------------------------------------------